\documentclass[12pt,a4paper]{article}
\usepackage{amssymb,amsfonts}
\usepackage{latexsym}
\usepackage{amsmath,amsthm}
\theoremstyle{plain}
\newtheorem{theorem}{Theorem}[section]
\newtheorem{corollary}{Corollary}
\newtheorem{lemma}{Lemma}
\newtheorem{proposition}{Proposition}
\newtheorem{example}{Example}
\theoremstyle{remark}
\newtheorem{remark}{Remark}[section]
\theoremstyle{definition}
\newtheorem{definition}{Definition}[section]

\newcommand{\R}{\mathbb R}
\newcommand{\E}{\mathbb R^n}
\newcommand{\f}{f:\Gamma\to\R}

\newcommand{\pr}{\prime}

\newcommand{\norm}[1]{\Vert#1\Vert}

\newtheorem*{KKT}{Karush-Kuhn-Tucker's Theorem}{\bf}{\it}
\newtheorem*{Gordan}{Gordan's Theorem of the Alternative}{\bf}{\it}
\title{Characterizations of Solution Sets of Fr\'echet Differentiable Problems with Quasiconvex Objective Function} 
\author{Vsevolod I. Ivanov\thanks{Email: vsevolodivanov@yahoo.com
\vspace{6pt} }
\\\vspace{6pt}{\em{\small Department of Mathematics, Technical University of Varna, 9010 Varna, Bulgaria} }
}
\date{\today}

\begin{document}
\maketitle

\begin{abstract}
In this paper, we study some problems with continuously differentiable quasiconvex  objective function. We prove that  exactly one of the following two alternatives holds: 
(I) the gradient of the objective function is different from zero over the solution set and the normalized gradient is constant over it;
(II) the gradient of the objective function is equal to zero over the solution set.
As a consequence, we obtain characterizations of the solution set of a quasiconvex continuously differentiable program, provided that one of the solutions is known. We also derive Lagrange multiplier characterizations of the solutions set of an inequality constrained problem with continuously differentiable objective function and differentiable constraints, which are all quasiconvex on some convex set, not necessarily open. We compare our results with the previous ones. Several examples are provided.

%\keywords{Global optimization \and Quasiconvex function} 



\end{abstract}

\section{Introduction}

In 1988, Mangasarian \cite{man88} obtained characterizations of the solution set of a convex program in terms of a known solution. 
 % Some more characterizations in terms of a known solution were derived by Jeyakumar and Yang \cite{jey95}, where the treated program were pseudolinear. 
%	Ivanov \cite{2001,2003} extended Jeyakumar and Yang's results to pseudoconvex programs. 
A lot of papers appeared later dealing with convex, pseudolinear, pseudoconvex, invex and other types of nonlinear programming problems. %Ivanov extended also Mangasarian's results concerning twice continuously differentiable convex programs to pseudoconvex ones (see Ref. \cite{jota2013}). 
%The Refs. \cite{jota2013,lon16,son14,wu08,zha13} concern pseudoconvex programs, Refs. \cite{din06,lal09,smi12} have treated pseudolinear ones. In Ref. \cite{jey02}, Jeyakumar, Lee, and Dinh obtained Lagrange multiplier characterizations for convex inequality constrained problems. Penot \cite{pen03} obtained some optimality conditions for quasiconvex programs in terms of several subdifferentials and invariance of the inward  subdifferential for essentially  quasiconvex programs.  In Ref. \cite{jota2010} were studied higher-order pseudoconvex programs. These are problems with a quasiconvex objective function, which is more general than pseudoconvex ones. Characterization of the approximate solution set of some problems were derived in \cite{liu15}.

Recently, Suzuki and Kuroiwa \cite{suz15} obtained generalizations of the Mangasarian's characterizations to characterizations of the solution set of an essentially quasiconvex program in terms of the Greenberg-Pierscalla subdifferential. 
It is well known that a real continuous function is essentially quasiconvex if and only if it is semistrictly quasiconvex \cite{avr88}. The Greenberg-Pierscalla subdifferential is defined for a quasiconvex  function. It does not necessarily includes the gradient of the function in the case when the last one is Fr\'echet differentiable. Therefore, it is important to derive simple characteriztions of solution set of an arbitrary  program with quasiconvex continuously differentiable objective function. Suzuki and Kuroiwa also derived an optimality condition and some more characterizations of the solution set of a quasiconvex program in terms of the Martinez-Legaz subdifferential (see \cite{suz17}), but it does not seem to be so simple.

In this paper, we continue the investigations initiated in the submission \cite{jmaa}. The same work
was submitted on February, 24th, 2010 to the new editor in chief of Journal of Optimization Theory and Applications, and it was rejected again by the associated editor on May 31th, 2011. In Ref. \cite{jota2013} are derived similar result, but the invex objective function were replaced by pseudoconvex one.  In the present paper, we obtain simple characteriztions of the solution set of a quasiconvex program with continuously differentiable objective function in terms of a given solution. We also extend some Lagrange multiplier characterizations of the solution set in terms of a given minimizer to quasiconvex programs. These results are known when the problem is convex or pseudoconvex.

%Lagrange multiplier characterizations were obtained in \cite{jey02} when the problem is convex, and generalized in \cite{zha13} to problems with pseudoconvex objective function and pseudoconvex constraints in terms of the Clarke generalized gradient.

\section{Characterizations of Solution Set of a Quasiconvex Program}
\label{s2}
\setcounter{theorem}{0}
\setcounter{definition}{0}
\setcounter{remark}{0}

Throughout this paper, $\E$ is the real $n$-dimensional Euclidean vector space,
$\Gamma\subseteq\E$ is an open convex set, $S\subseteq \Gamma$ is a nonempty convex subset of $\Gamma$, and $f$ is a Fr\'echet differentiable function, defined on $\Gamma$. The main purpose of this section is to obtain characterizations of solution set of the nonlinear programming problem:

\bigskip
Minimize $f(x)\quad$ subject to $\quad x\in S$. \hfill  (P)
\bigskip

\par
Denote the transpose of the matrix $A$ by $A^T$, the scalar product of the vectors $a\in\R^n$ and $b\in\R^n$ by $a^T\,b$, the solution set  $\arg\min\;\{f(x)\mid
x\in S\}$ of (P)  by $\bar S$, and let $\bar S$  be nonempty. Suppose that $\bar x$ is a known point from $\bar S$.

Recall the following well-known definition \cite{man69}:

\begin{definition}
A function $\f$, defined on a convex set $\Gamma\subseteq\E$, is called {\it quasiconvex} on $\Gamma$ iff the following condition holds for all $x$, $y\in\Gamma$, $t\in[0,1]$:
\[
f[x+t(y-x)]\le\max\{f(x),f(y)\}.
\]
\end{definition}

We begin with some preliminary results:

\begin{lemma}[\cite{man69}]\label{lema1}
A function $\f$ is quasiconvex on the convex set $\Gamma\subseteq\E$ if and only if its lower level sets $L(f;\alpha):=\{x\in\Gamma\mid f(x)\le\alpha\}$ are convex for all real numbers $\alpha$.
\end{lemma}

\begin{lemma}[\cite{arr61}]\label{lema2}
Let $f:\Gamma\to\R$ be a quasiconvex differentiable function on the open convex set $\Gamma$. Suppose that $x\in\Gamma$, $y\in\Gamma$, and
$f(y)\le f(x)$. Then $\nabla f(x)^T(y-x)\le 0$  
\end{lemma}

\begin{Gordan}{\rm\cite{gor1873,man69}}
For each given matrix $A$, either the system
\[
A\,x>0
\]
has a solution $x$, or the system
\[
A^T\,y=0,\; y\ge 0,\; y\ne 0
\]
has a solution $y$, but never both.
\end{Gordan}

\begin{lemma}\label{lema3}
Let $\Gamma\subseteq\E$ be an open convex set, $S\subseteq\Gamma$ be a convex set and $x$, $y\in\bar S$ with $\nabla f(x)\ne 0$, $\nabla f(y)\ne 0$. Suppose that $f$ is quasiconvex on $\Gamma$. Then the following implication holds:
\[
\nabla f(x)^Td<0\quad\Rightarrow\quad\nabla f(y)^Td\le0.
\]
\end{lemma}
 \begin{proof}%\smartqed
Suppose the contrary that there exists $d\in\E$ with  
\[
\nabla f(x)^Td<0\quad{\rm and}\quad \nabla f(y)^T(-d)<0.
\]
 It follows from here that there exists $\tau>0$ with
\begin{equation}\label{16}
f(u)<f(x)=f(y),\quad f(v)<f(y)=f(x),
\end{equation}
where $u=x+\tau d\in\Gamma$ and $v=y-\tau d\in\Gamma$. Let $z=(u+v)/2$. Therefore $z=(x+y)/2$.  According to the quasiconvexity of $f$, we obtain by Lemma \ref{lema1}  that $\bar S$ is convex and $f(z)=f(x)=f(y)$. By quasiconvexity, we conclude from (\ref{16}) that
\begin{equation}\label{17}
f(z)\le\max\{f(u),f(v)\}<f(x)=f(z),
\end{equation}
which is impossible.
\end{proof}

\begin{lemma}\label{lema4}
Let $a$, $b\in\R^n$, $a\ne 0$, $b\ne 0$. Suppose that
\begin{equation}\label{4}
a^T\,d<0,\; d\in\R^n \quad\Rightarrow\quad b^T\, d\le 0.
\end{equation}
Then there exists $p>0$ such that $b=p\, a$.
\end{lemma}
\begin{proof}   %\smartqed
Implication (\ref{4}) is equivalent to the claim that the system
\[
b^T\, d>0,\quad   (-a)^T\,d>0
\]
has no a solution $d$. It follows from Gordan's Theorem that there exist real numbers $p_1$ and $p_2$ such that
\[
p_1\, b-p_2\, a=0,\quad p_1\ge 0,\; p_2\ge 0,\; (p_1,p_2)\ne (0,0).
\]
Without loss of generality $p_1>0$. Let $p=p_2/p_1$. The number $p$ is strictly positive, and it satisfies the equation $b=p\, a$, because $b\ne 0$. %The proof is complete.
\end{proof}

\begin{lemma}\label{lema5}
Let $\Gamma\subseteq\E$ be an open convex set, $S\subseteq\Gamma$ be a convex set and $f:\Gamma\to\R$ is a Fr\'echet differentiable quasiconvex function. Then the normalized gradient is constant over the set $\{x\in\bar S\mid\nabla f(x)\ne 0\}$, provided that this set is nonempty.
\end{lemma}
\begin{proof}
Let $x$ and $y$ be arbitrary points from the set $\bar S$ with $\nabla f(x)\ne 0$, $\nabla f(y)\ne 0$. It follows from Lemmas \ref{lema3} and \ref{lema4} that there exists $p(x,y)>0$ with $\nabla f(y)=p\nabla f(x)$.  We obtain from here that
\[
\nabla f(y)/\norm{\nabla f(y)}=p\nabla f(x)/\norm{p\nabla f(x)}=\nabla f(x)/\norm{\nabla f(x)}.
\]
Then the claim follows immediately.
\end{proof}	

\begin{lemma}\label{lema6}
Let $\Gamma\subseteq\E$ be an open convex set, $S\subseteq\Gamma$ be a convex one. Suppose that $f:\Gamma\to\R$ is a continuously differentiable quasiconvex function. Then exactly one of the alternatives holds:

I) $\nabla f(x)\ne 0$ for all $x\in\bar S$ and the normalized gradient $\nabla f(x)/\norm{\nabla f(x)}$ is constant over the solution set $\bar S$;

II) $\nabla f(x)=0$ for all $x\in\bar S$.
\end{lemma}
\begin{proof}
Let us take an arbitrary point $x$ from the solution set $\bar S$. Two cases are possible:

First case: $\nabla f(x)\ne 0$. Let $y$ be another arbitrary point from $\bar S$. We prove that $\nabla f(y)\ne 0$. Suppose the contrary that $\nabla f(y)=0$. Consider the sets:
\[
A:=\{t\in[0,1]\mid\nabla f[x+t(y-x)]=0\},\quad B:=\{t\in[0,1]\mid\nabla f[x+t(y-x)]\ne 0\}.
\]
We have $x+t(y-x)\in\bar S$, because $\bar S$ is convex by Lemma \ref{lema2}. Consider the standard topology on the interval $[0,1]$. The open sets in it are the open intervals $(a,b)$ such that $0<a<b<1$, the intervals $(a,1]$ such that $0<a<1$, the intervals $[0,b)$ such that $0<b<1$and their unions. 

We prove that $A$ is closed. Let us take a sequence $\{t_n\}$, where $t_n\in A$ and $t_n$ approaches to $t_0$. By the continuous differentiability of $f$ we have
\[
\nabla f[x+t_0(y-x)]=\lim_{n\to\infty}\nabla f[x+t_n(y-x)]=\lim_{n\to\infty}\, 0=0.
\]
Therefore $t_0\in A$ and $A$ is a closed set.

We prove that $B$ is a closed set. Let us take a sequence $\{t_n\}$, where $t_n\in B$ and $t_n$ approaches to $t_0$. By the continuous differentiability of $f$ we have
\[
\frac{\nabla f[x+t_0(y-x)]}{\norm{\nabla f[x+t_0(y-x)]}}=\lim_{n\to\infty}\frac{\nabla f[x+t_n(y-x)]}{\norm{\nabla f[x+t_n(y-x)]}},
\]
because by  Lemma \ref{lema5} the normalized gradient is a constant vector, different from $0$. Therefore $t_0\in B$ and $B$ is closed.

Since the union of $A$ and $B$ is the whole interval $[0,1]$, then $A=[0,1]$, $B=\emptyset$, or $B=[0,1]$, $A=\emptyset$. Both cases are impossible, because $\nabla f(x)\ne 0$ and $\nabla f(y)=0$. This fact is contrary to the assumption $\nabla f(y)=0$. Therefore $\nabla f(y)\ne 0$. It follows from here that $\nabla f(y)\ne 0$ for all $y\in\bar S$. By Lemma \ref{lema5} the normalized gradient of $f$ is constant over $\bar S$.

Second case: $\nabla f(x)=0$. It follows from the proof of the first case the impossibility of the assumption $\nabla f(y)\ne 0$ for arbitrary $y\in \bar S$. Therefore $\nabla f(y)=0$ for all $y\in\bar S$.
\end{proof}

Let $\bar x\in\bar S$ be a given point. Denote
\[\begin{array}{l}
\hat S_1:=\{x\in S\mid\nabla f(\bar x)^T(x-\bar x)=0,\;  \frac{\nabla f(x)}{\norm{\nabla f(x)}}=\frac{\nabla f(\bar x)}{\norm{\nabla f(\bar x)}},\;\nabla f(x)\ne 0\}; \\
%\]
%\[
\hat S_2:=\{x\in S\mid\nabla f(\bar x)^T(x-\bar x)\le 0,\;\frac{\nabla f(x)}{\norm{\nabla f(x)}}=\frac{\nabla f(\bar x)}{\norm{\nabla f(\bar x)}},\;\nabla f(x)\ne 0\}.
\end{array}\]

\begin{theorem}\label{th1}
Let $\Gamma\subseteq\E$ be an open convex set, $S\subseteq\Gamma$ be a convex one and $f:\Gamma\to\R$ be a continuously differentiable quasiconvex function. Suppose that $\bar x\in\bar S$ is a known solution of {\rm (P)} and $\nabla f(\bar x)\ne 0$. Then
\[
\bar S=\hat S_1=\hat S_2.
\]
\end{theorem}
 \begin{proof}  %\smartqed
We prove that $\bar S\subseteq\hat S_1$. Suppose that $x\in\bar S$. We prove that $x\in\hat S_1$. By Lemma \ref{lema1}, $\bar S$ is convex and we have $\bar x+t(x-\bar x)\in\bar S$ for all $t\in[0,1]$. Therefore 
\[
f[\bar x+t(x-\bar x)]=f(\bar x),\quad\forall t\in[0,1],
\]
because the solution set $\bar S$ is convex. It follows from here that $\nabla f(\bar x)^T(x-\bar x)=0$.
According to Theorem \ref{lema6} $\nabla f(x)\ne 0$ and $x\in\hat S_1$.
 
It is trivial that $\hat S_1\subseteq\hat S_2$ 

We prove that $\hat S_2\subseteq\bar S$. Let $x\in\hat S_2$. We prove that $x\in\bar S$. Suppose the contrary. Therefore $f(\bar x)<f(x)$. According to the convexity of $S$ we have $\bar x+t(x-\bar x)\in S$. By the assumption $\bar x\in\bar S$, we obtain that 
$f[\bar x+t(x-\bar x)]\ge f(\bar x)$. Therefore
 $\nabla f(\bar x)^T(x-\bar x)\ge 0$. It follows from  $x\in\hat S_1$ that $\nabla f(\bar x)^T(x-\bar x)=0$. Since the function $f$ is continuous and $f(\bar x)<f(x)$, there exists a number $\delta>0$ such that $f[\bar x+\delta\nabla f(\bar x)]<f(x)$. Then, we conclude from Lemma  \ref{lema2} that
\[
\nabla f(x)^T[\bar x+\delta\nabla f(\bar x)-x]\le 0.
\]
%According to Lemma \ref{lema6} $\nabla f(x)\ne 0$.
By the equality $\nabla f(x)/\norm{\nabla f(x)}=\nabla f(\bar x)/\norm{\nabla f(\bar x)}$ we obtain that  
\[
\nabla f(\bar x)^T[\bar x+\delta\nabla f(\bar x)-x]\le 0.
\]
It follows from $\nabla f(\bar x)^T(x-\bar x)=0$ that $\norm{\nabla f(\bar x)}^2\le 0$, which contradicts the hypothesis $\nabla f(\bar x)\ne 0$.
\end{proof}

%We can see immediately that Theorem \ref{th1} improves and generalizes this result to quasiconvex objective functions. We can say that a quasiconvex function, whose gradient does not vanish on the fesible set is really pseudoconvex on $S$.

\begin{example}\label{ex2}
Consider the problem (P) such that function $f$ is the function of two variables $f(x_1,x_2)=x_2/x_1$ and $S$ is the rectangle
\[
S=\{x=(x_1,x_2)\mid 1\le x_1\le 2,\; 0\le x_2\le x_1\}.
\]
 The objective function is quasiconvex on the  set $\Gamma=\{x=(x_1,x_2)\mid x_1>0\}$, because it satisfy the implication
\[
y\in\Gamma,\; x\in\Gamma,\; f(y)\le f(x)\quad\Rightarrow\quad \nabla f(x)^T(y-x)\le 0.
\]
The solution set is $\bar S=\{(x_1,x_2)\mid 1\le x_1\le 2,\; x_2=0\}$. Let $\bar x=(1,0)$ be a known solution. We can find $\bar S$ applying Theorem \ref{th1}, because $\hat S_1=\{(x_1,x_2)\mid 1\le x_1\le 2,\; x_2=0\}$.
\end{example}

The  objective function of the previous problem is really pseudoconvex. In the next example $f$ is not pseudoconvex.

\begin{example}\label{ex3}
Consider the problem (P) such that function $f$ is the function of two variables $f(x_1,x_2)=x_1^3$ and 
\[
S:=\{x=(x_1,x_2)\mid x_1\ge -1,\; -\infty<x_2<+\infty\}.
\]
$f$ is quasiconvex, but not pseudoconvex. Let $\bar x=(-1,0)$ be a known solution. It is easy to see that $\hat S=\{(x_1,x_2)\mid x_1=-1\}$. By Theorem \ref{th1} this is the solution set.
\end{example}

\begin{example}\label{ex4}
Consider the problem (P) such that function $f$ is the function of two variables $f(x_1,x_2)=-x_1-x_2+\sqrt{(x_1-x_2)^2+4}$ and $S$ is the disk
\[
S=\{x=(x_1,x_2)\mid x_1^2+x_2^2\le 2\}.
\]
The objective function is quasiconvex, because its lower level set are convex. Really, for every $c\in\R$ the inequality $f(x)\le 2c$ is equivalent to
\[
(x_1+c)(x_2+c)\ge 1,\; x_1+x_2+2c\ge 0.
\]
Let $\bar x=(1,1)$ be a known solution. It is easy to see that $\hat S$ consists of the only point $(1,1)$. By Theorem \ref{th1}, the problem has no more solutions.
\end{example}

%Theorem \ref{th1} together with Theorem \ref{th2} below are  generalization of Theorem 1 in the paper by Mangasarian \cite{man88}, where the function is twice continuously differentiable and convex. It is also a generalization of Theorem 2.1 in the paper by Ivanov \cite{jota2013}, where the objective function is pseudoconvex.

%Consider the following result due to Jeyakumar and Yang \cite[Proposition 3.1]{jey95}:
%\begin{proposition}
%Let $f$ be continuously differentiable and pseudolinear on an open convex set containing the convex set $S\subseteq\R^n$; let 
%$\bar x\in\Bar S$. Suppose that $\nabla f(x)\ne 0$ on $S$. Then
%\[
%\bar S\subseteq\{x\in S\mid\nabla f(x)/\norm{\nabla f(x)}=\nabla f(\bar x)/\norm{\nabla f(\bar x)}\}\subseteq\{x\in S\mid\nabla f(x)^T(\bar x-x)\le 0\}.
%\]
%\end{proposition}

%We can see immediately that Theorem \ref{th1} improves and generalizes this result to quasiconvex objective functions. We can say that a quasiconvex function, whose gradient does not vanish on the fesible set is really pseudoconvex on $S$.

\begin{definition}[\cite{man69}]
Let $\Gamma\subseteq\R^n$ be an open set. A differentiable function $f:\Gamma\to\R$ is called pseudoconvex on $\Gamma$ if the following implication holds:
\[
x\in\Gamma,\; y\in\Gamma,\; f(y)<f(x)\quad\Rightarrow\quad \nabla f(x)^T(y-x)<0.
\]
If this implication is satisfied at the point $x$ only for every $y\in\Gamma$, then the function is said to be pseudoconvex at $x$.  
\end{definition}
Every differentiable function, pseudoconvex on some open convex set, is quasiconvex \cite{man69}.

\smallskip
Consider the following set:
\[
\tilde S:=\{x\in S\mid\nabla f(x)=0\}.
\]

\begin{theorem}\label{th2}
Let $\Gamma\subseteq\E$ be an open convex set, $S\subseteq\Gamma$ be a convex one and $f:\Gamma\to\R$ be a continuously differentiable quasiconvex function. Suppose that $\bar x\in\bar S$ is a known solution of {\rm (P)} and $\nabla f(\bar x)=0$. Then $\nabla f(x)=0$ for all $x\in\bar S$. 

Suppose additionally that $f$ is pseudoconvex at every point $x\in S$ such that $x\notin\bar S$.  Then $\bar S=\tilde S$.
\end{theorem}
\begin{proof}   %\smartqed
%The theorem follows directly from Lemma \ref{lema6}.
The claim $\bar S\subseteq\tilde S$ follows directly from Lemma \ref{lema6}. It follows from this lemma that the second alternative holds in our case. Therefore, $\nabla f(x)=0$ for all $x\in\bar S$.

Suppose that $f$ is pseudoconvex on $S\setminus\bar S$. We prove that $\tilde S\subseteq\bar S$. Let $x$ be an arbitrary point from $\tilde  S$. Assume the contrary that $x\notin\bar S$. It follows from here that $f(\bar x)<f(x)$. By pseudoconvexity, we obtain that 
$\nabla f(x)^T(\bar x-x)<0$. This inequality contradicts the condition $\nabla f(x)=0$, which follows from $x\in\tilde S$.

\end{proof}

\begin{example}\label{ex1}
Consider the function $f:\R^2\to\R$, defined as follows:
\[
f(x)=\left\{\begin{array}{ll}
x_1^2+x_2^2, & x_1\ge 0,\; x_2\ge 0, \\
x_2^2, & x_1\le 0,\; x_2\ge 0, \\ 
-x_1^2 x_2^2, & x_1\le 0,\; x_2\le 0, \\ 
x_1^2, & x_1\ge 0,\; x_2\le 0
\end{array}\right.
\]
and the problem (P) such that
%\[
%\textrm{Minimize}\quad f(x)\quad\textrm{subject to}\quad x\in S,
%\]
%where 
$S:=\{x=(x_1,x_2)\in\R^2\mid x_1\ge 0\}$. The objective function is quasiconvex, because its lower level sets are convex. Really, it is convex over $S$. On the other hand Theorem 1 from \cite{man88} cannot be applied, because $S$ is not open and there is no an open convex set such that $f$ is convex over it. It is easy to see that $f\in$C$\,^1$  and it is not pseudoconvex. Therefore, the theory concerning pseudoconvex programs also cannot be applied. The solution set is $\{x=(x_1,x_2)\in\R^2\mid x_1=0,\; x_2\le 0\}$. It is easy to see that the function is pseudoconvex at every $x\in\{x=(x_1,x_2)\mid x_1>0\}\cup \{x=(x_1,x_2)\mid x_1=0,\; x_2>0\}$. 
Let $\bar x=(0,0)$ be a known solution. It follows from theorem \ref{th2} that 
\[
\bar S=\{x\in S\mid\nabla f(x)=0\}=\{x\in\R^2\mid x_1=0,\; x_2\le 0\}.
\]
%Тази функция е псевдоизпъкнала от четвърти ред и диференцируема навсякъде. Това се доказва подобно на функцията от пример \ref{6.ex4}. Когато доказваме импликацията $f(y)<f(x)$ означава, че $\nabla f(x)(y-x)<0$ разглеждаме същите случаи, както и в предния пример. Различен е единствено случаят, когато $x$ и $y$ са точки от трети  квадрант. От неравенството $f(y)<f(x)$ следва, че $x_1 x_2<y_1 y_2$. Тъй като $y_1\ne 0$ оттук получаваме, че $y_2<x_1 x _2/y_1$. Като вземем предвид, че $\nabla f(x)(y-x)=2 x_1 x_2(2x_1 x_2-x_1 y_2-x_2 y_1)$ от неравенството $y_2<x_1 x _2/y_1$ лесно се вижда, че $\nabla f(x)(y-x)<0$. Ако $x=(0,0)$ и $y$ е точка от трети квадрант, то $f^\pr(x,y-x)=f^{\pr\pr}(x,y-x)=f^{\pr\pr\pr}(x,y-x)=0$ и $f^{\pr v}(x,y-x)<0$. Ако $x_1=0$, $x_2<0$ или $x_1<0$, $x_2=0$, то $\nabla f(x)(y-x)=0$ и $f^{\pr\pr}(x,y-x)<0$.
\end{example}

Consider the following sets:
\[\begin{array}{l}
S_1:=\{x\in S\mid \nabla f(x)^T(\bar x-x)=0,\; \nabla f(x)\ne 0\}, \\
%\]
%\[
S_2:=\{x\in S\mid \nabla f(x)^T(\bar x-x)\ge 0,\; \nabla f(x)\ne 0\}, \\
%\]
%\[
S_3:=\{x\in S\mid \nabla f(x)^T(\bar x-x)=\nabla f(\bar x)^T(x-\bar x),\; \nabla f(x)\ne 0\}, \\
%\]
%\[
S_4:=\{x\in S\mid \nabla f(x)^T(\bar x-x)\ge\nabla f(\bar x)^T(x-\bar x),\; \nabla f(x)\ne 0\}, \\
%\]
%\[
S_5:=\{x\in S\mid\nabla f(x)^T(\bar x-x)=\nabla f(\bar x)^T(x-\bar x)=0,\; \nabla f(x)\ne 0\}. 
\end{array}\]

\begin{theorem}\label{th3}
 Let $\Gamma\subseteq\E$ be an open convex set, $S\subseteq\Gamma$ be an arbitrary convex one and the function $f$ be defined on $\Gamma$.
Suppose $\bar x\in\bar S$ is a known solution of {\rm (P)} such that $\nabla f(\bar x)\ne 0$.
If the function $f$ is continuously differentiable and  quasiconvex on $\Gamma$, then
\[
\bar  S=S_1=S_2=S_3=S_4=S_5.
\]
\end{theorem}
 \begin{proof} %\smartqed
It is obvious that $S_5\subseteq S_1\subseteq S_2$ and $S_5\subseteq S_3\subseteq S_4$.

We prove that  $\bar  S\subseteq S_5$. Suppose that $x\in\bar S$. Therefore $f(x)=f(\bar x)$. By quasiconvexity the solution set $\bar S$ is convex and
\[
f[\bar x+t(x-\bar x)]=f(\bar x)\quad\textrm{for all}\quad t\in [0,1].
\]
Therefore $\nabla f(\bar x)^T(x-\bar x)=0$. We can prove using similar arguments,  interchanging $\bar x$ and $x$, that $\nabla f(x)^T(\bar x-x)=0$. The claim that $\nabla f(x)\ne 0$ follows from Lemma \ref{lema6}, because by the hypothesis $\nabla f(\bar x)\ne 0$.
We conclude from here that $x\in S_5$.

We prove that $S_4\subseteq S_2$. Let $x\in S_4$. Therefore 
\begin{equation}\label{6}
\nabla f(x)^T(\bar x-x)\ge\nabla f(\bar x)^T(x-\bar x).
\end{equation}
Since $S$ is convex, $x$, $\bar x\in S$ and $\bar x\in\bar S$, then $f[\bar x+t(x-\bar x)]\ge f(\bar x)$. Therefore 
\[
\nabla f(\bar x)^T(x-\bar x)\ge 0,
\]
which implies by (\ref{6}) that  $x\in S_2$.
%$\nabla(x)(\bar x-x)\ge 0$.

At last, we prove that $S_2\subseteq\bar S$. Let $x\in S_2$. Assume the contrary that $x\notin\bar S$. Hence $f(\bar x)<f(x)$. By Lemma \ref{lema2}, we obtain that $\nabla f(x)^T(\bar x-x)\le 0$.  Taking into account that $x\in S_2$ we have $\nabla f(x)^T(\bar x-x)=0$.  Then, by the continuity of $f$, it follows from $f(\bar x)<f(x)$ that there exists $\delta>0$ with
$f[\bar x+\delta\nabla f(x)]<f(x)$. Then, it follows from Lemma \ref{lema2} that 
\[
\nabla f(x)^T(\bar x+\delta\nabla f(x)-x)\le 0.
\]
By  $\nabla f(x)^T(\bar x-x)=0$, we conclude from here that $\nabla f(x)=0$, which is a contradiction to $x\in S_2$. 
\end{proof}

%Theorem \ref{th3} is a generalization of the main results  in the paper \cite{jey95} by Jeyakumar and Yang (see Theorem 3.1, Theorem 3.2 and Corollary 3.1), where the objective function is differentiable and pseudolinear. Characterizations, which are similar to Theorem \ref{th3} were obtained in the papers by Ivanov \cite[Theorem 4.1]{2001}, \cite[Theorem 3.1]{2003}, where the objective function is nondifferentiable and pseudoconvex .

\section{Lagrange multiplier characterizations of the solution set}
Consider the problem with inequality constraints

\bigskip
Minimize $f(x)\quad$ subject to $x\in X,\quad g_i(x)\le 0$, $i=1,2,...,m$,
\hfill {\rm (PI)}
\bigskip

\noindent
where $f:\Gamma\to\R$ and $g_i:\Gamma\to\R$ are defined on some open set $\Gamma\subseteq\R^n$, $X$ is a convex subset of $\Gamma$, not necessarily open.

Denote  by
\[
I(x):=\{i\in\{1,2,...,m\}\mid g_i(x)=0\}
\]
the index set of the active constraints at the feasible point $x$. Let 
\[
S:=\{x\in X\mid g_i(x)\le 0,\; i=1,2,...,m\}
\]
be the feasible set.

Suppose that $C\subseteq\R^n$ is a cone. Then the cone $C^*:=\{x\in\R^n\mid c^T\, x\le 0\}$ is said to be the negative polar cone of $C$.  Let $T_X(x)$ be the tangent cone of the set $X$ at the point $x$. Then its negative polar cone is called the normal cone $N_X(x)$.

\begin{definition}\cite{ber03}
It is said that the constraint functions satisfy generalized Mangasari\-an-Fromovitz constraint qualification (in short, GMFCQ) at the point $\bar x$ iff there is a direction $y\in (N_X(\bar x))^*$ such that $\nabla g_i(\bar x)^Ty<0$ for all $i\in I(\bar x)$.
%the set $\Gamma$ is open and the vectors $\nabla g_i(\bar x)$, $i\in I(\bar x)$ are linearly independent.
\end{definition}
The following necessary optimality conditions of Karush-Kuhn-Tucker type (for short, KKT conditions) are consequence of Proposition 2.2.1, Definition 2.4.1 and Proposition 2.4.1 in \cite{ber03}:
%well known (see, for example, \cite{man69}):

\begin{KKT}\label{KKT}
Let $\bar x$ be a local minimizer of the problem {\rm (PI)}. Suppose that $f$, $g_i$, $i\in I(\bar x)$ are Fr\'echet differentiable on 
$\Gamma\subseteq\R^n$ at $\bar x$, $g_i$, $i\notin I(\bar x)$ are continuous at $\bar x$, the set $X$ is convex. Suppose additionally that GMFCQ holds at $\bar x$.  Then there exists a Lagrange multiplier $\lambda\in\R^m$,
$\lambda=(\lambda_1,\lambda_2,...,\lambda_m)$, $\lambda_i\ge 0$, $i=1,2,\dots,m$ such that
\[
[\nabla f(\bar x)+\sum_{i\in I(\bar x)}\lambda_i\nabla g_i(\bar x)]^T(x-\bar x)\ge 0,\;\forall\, x\in X,\quad \lambda_i g_i(\bar x)=0\; \forall\; i=1,2,...,m.
\]
\end{KKT}

Denote  by $\tilde I(\bar x,\lambda)$ the following index set
\[
\tilde I(\bar x,\lambda):=\{i\in\{1,2,...,m\}\mid g_i(\bar x)=0,\;\lambda_i>0\}
\]
and the set
\[
X_1(\lambda):=\{x\in X\mid g_i(x)=0\;\forall i\in\tilde I(\bar x,\lambda),\;\; g_i(x)\le 0\;\forall i\in\{1,2,\dots,m\}\setminus\tilde I(\bar x,\lambda)\}.
\]

\begin{lemma}\label{lema7}
Let the functions $f$ and $g$ be differentiable and quasiconvex, $\bar x\in\bar S$ be a solution. Suppose that the set $X$ is convex,  GMFCQ is satisfied at $\bar x$ and KKT optimality conditions are satisfied at $\bar x$ with a multiplier $\lambda$.
Then $\bar S\subseteq X_1(\lambda)$ and the Lagrangian function $L=f(\cdot)+\sum_{i\in I(\bar x)}\lambda_i g_i(\cdot)$ is constant over $\bar S$.
\end{lemma}
\begin{proof}   %\smartqed
We prove that $\bar S\subseteq X_1(\lambda)$.
 Let $x$ be an arbitrary element of $\bar S$. We prove that $x\in X_1(\lambda)$.  It is enough to show that $g_i(x)=0$ for all $i\in\tilde I(\bar x,\lambda)$, because $x$ is a feasible point. Suppose the contrary that there exists $j\in\tilde I(\bar x,\lambda)$ such that $g_j(x)<0$. Taking into account that $g_j(\bar x)=0$, then we have $g_j(x)<g_j(\bar x)$. There exists $\delta>0$ such that $g_j[x+\delta\nabla g_j(\bar x)]<g_j(\bar x)$. Then it follows from Lemma \ref{lema2} that 
\begin{equation}\label{3}
\nabla g_j(\bar x)^T[x+\delta\nabla g_j(\bar x)-\bar x]\le 0.
\end{equation}
By quasiconvexity of $f$ and $f(x)=f(\bar x)$ it follows from Lemma \ref{lema2} that $\nabla f(\bar x)^T(x-\bar x)\le 0$. Since $g_i$ are also quasiconvex and $g_i(x)\le 0=g_i(\bar x)$ for all $i\in I(\bar x)$, by Lemma \ref{lema2} we obtain that $\nabla g_i(\bar x)^T(x-\bar x)\le 0$ for all $i\in I(\bar x)$. By KKT conditions the following equations are satisfied:
\begin{equation}\label{1}
\nabla f(\bar x)^T(x-\bar x)=0,\quad \lambda_i\nabla g_i(\bar x)^T(x-\bar x)=0\; \forall\; i\in I(\bar x).
\end{equation}
 By (\ref{1}), we obtain that $\nabla g_j(\bar x)^T(x-\bar x)=0$.     Therefore, by (\ref{3}) we have
 $\norm{\nabla g_j(\bar x)}^2\le 0$. On the other hand, it follows from GMFCQ that $\nabla g_j(\bar x)\ne 0$, which contradicts the indirect conclusion that $\nabla g_j(\bar x)=0$. Hence $x\in  X_1(\lambda)$. 

The claim that $L(x)=L(\bar x)$ follows immediately from here according to the equality $f(x)=f(\bar x)$.
;\end{proof}

Consider the sets
\[\begin{array}{l}
\hat S_1^\pr(\lambda):=\{x\in X_1(\lambda)\mid\nabla f(\bar x)^T(x-\bar x)=0,\;\;\frac{\nabla f(x)}{\norm{\nabla f(x)}}=\frac{\nabla f(\bar x)}{\norm{\nabla f(\bar x)}},\;\nabla f(x)\ne 0\}; \\
%\]
%\[
\hat S_2^\pr(\lambda):=\{x\in X_1(\lambda)\mid\nabla f(\bar x)^T(x-\bar x)\le 0,\;\;\frac{\nabla f(x)}{\norm{\nabla f(x)}}=\frac{\nabla f(\bar x)}{\norm{\nabla f(\bar x)}},\;\nabla f(x)\ne 0\}.
\end{array}\]

\begin{theorem}\label{th4}
Let the set $X$ be  convex, the function $f$ be continuously differentiable and quasiconvex,  $g_i$, $i\in I(\bar x)$ be differentiable and quasiconvex,  $g_i$, $i\notin I(\bar x)$ be continuous at $\bar x$. Suppose that  $\bar x\in\bar S$, $\nabla f(\bar x)\ne 0$, GMFCQ is satisfied, the Lagrange multipliers are known and fixed. Then
\[
\bar S=\hat S_1^\pr(\lambda)=\hat S_2^\pr(\lambda).
\]
\end{theorem}
\begin{proof}   %\smartqed
It is obvious that $\hat S_1^\pr(\lambda)=X_1(\lambda)\cap\hat S_1$ and $\hat S_2^\pr(\lambda)=X_1(\lambda)\cap\hat S_2$. Then the claim follows from Lemma \ref{lema7} and the equality $\bar S=\hat S_1=\hat S_2$, which is a part of Theorem \ref{th1}.
\end{proof}

\begin{remark}
In Theorem \ref{th4}, we suppose that $\nabla f(\bar x)\ne 0$, but if $\nabla f(\bar x)=0$, then  $\nabla f(x)=0$ for all $x\in\bar S$ according to Theorem \ref{th2}.
\end{remark}

Consider the sets
\[\begin{array}{l}
S_1^\pr(\lambda):=\{x\in X_1(\lambda)\mid\nabla f(x)^T(\bar x-x)=0,\;\; \nabla f(x)\ne 0\}; \\
%\]
%\[
S_2^\pr(\lambda):=\{x\in X_1(\lambda)\mid\nabla f(x)^T(\bar x-x)\ge 0,\;\; \nabla f(x)\ne 0\}; \\
%\]
%\[
S_3^\pr(\lambda):=\{x\in X_1(\lambda)\mid\nabla f(x)^T(\bar x-x)=\nabla f(\bar x)^T(x-\bar x),\;\; \nabla f(x)\ne 0\}; \\
%\]
%\[
S_4^\pr(\lambda):=\{x\in X_1(\lambda)\mid\nabla f(x)^T(\bar x-x)\ge \nabla f(\bar x)^T(x-\bar x),\;\; \nabla f(x)\ne 0\}; \\
%\]
%\[
S_5^\pr(\lambda):=\{x\in X_1(\lambda)\mid\nabla f(x)^T(\bar x-x)=\nabla f(\bar x)^T(x-\bar x)=0,\;\; \nabla f(x)\ne 0\}.
\end{array}\]

\begin{theorem}\label{th5}
Let the set $X$ be convex, the function $f$ be continuously differentiable and quasiconvex,  $g_i$, $i\in I(\bar x)$ be differentiable and quasiconvex, $g_i$, $i\notin I(\bar x)$ be continuous at $\bar x$. Suppose that $\bar x\in\bar S$, $\nabla f(\bar x)\ne 0$, GMFCQ is satisfied, the Lagrange multipliers are known and fixed. Then
\[
\bar S=S_1(\lambda)^\pr=S_2^\pr(\lambda)=S_3^\pr(\lambda)=S_4^\pr(\lambda)=S_5^\pr(\lambda).
\]
\end{theorem}
\begin{proof}   %\smartqed
It is obvious that $S_i^\pr(\lambda)=X_1(\lambda)\cap S_i$, $i=1,2,3,4,5$.  Then the claim follows from Lemma \ref{lema7} and the equalities $\bar S=S_i(\lambda)$, $i=1,2,3,4,5$ which are the statements of Theorem \ref{th3}.
\end{proof}

%Similar results to Lemma \ref{lema7} and Theorem \ref{th5} appear in \cite{jey02} in the case when the problem is convex, not necessarily differentiable (see Theorem 2.1 and 2.2). The reader can also compare Lemma \ref{lema7} and Theorem \ref{th5} with Theorems 3.1 and 4.1 in \cite{zha13}, where the problem contains locally Lipschitz pseudoconvex objective function and locally Lipschitz pseudoconvex constraints with respect to the Clarke generalized gradient, but the set $X$ is open and convex.

In the next result, we suppose that $X$ is an open set. Therefore KKT conditions reduce to the following one:
\begin{equation}\label{7}
\nabla f(\bar x)+\sum_{i\in I(\bar x)}\lambda_i\nabla g_i(\bar x)=0,\quad \lambda_i g_i(\bar x)=0\; \forall\; i=1,2,...,m.
\end{equation}

Generalized Mangasarian-Fromovitz constraint qualification reduces to Mangasarian-Fromovitz constraint qualification (in short, MFCQ): there is a direction $y\in\R^n$ such that $\nabla g_i(\bar x)^Ty<0$ for all $i\in I(\bar x)$.

In this case, we consider the following sets: 
\[\begin{array}{l}
\hat S_1^{\pr\pr}(\lambda):=\{x\in X_1(\lambda)\mid\nabla g_i(\bar x)^T(x-\bar x)=0,\; i\in\tilde I(\bar x,\lambda),\;
\frac{\nabla f(x)}{\norm{\nabla f(x)}}=\frac{\nabla f(\bar x)}{\norm{\nabla f(\bar x)}},\;\nabla f(x)\ne 0\}; \\
%\exists p(x)>0: \nabla f(x)=p(x)\nabla f(\bar x)\};
%\]
%\[
\hat S_2^{\pr\pr}(\lambda):=\{x\in X_1(\lambda)\mid\nabla g_i(\bar x)^T(x-\bar x)\ge 0,\; i\in\tilde I(\bar x,\lambda),\;\frac{\nabla f(x)}{\norm{\nabla f(x)}}=\frac{\nabla f(\bar x)}{\norm{\nabla f(\bar x)}},\;\nabla f(x)\ne 0\}.
%\exists p(x)>0: \nabla f(x)=p(x)\nabla f(\bar x)\};
\end{array}\]

\begin{theorem}\label{th6}
Let the set $X$ be open and convex, the function $f$ be continuously differentiable and quasiconvex,  $g_i$, $i\in I(\bar x)$ be differentiable and quasiconvex, $g_i$, $i\notin I(\bar x)$ be continuous at $\bar x$. Suppose that  $\bar x\in\bar S$, $\nabla f(\bar x)\ne 0$, MFCQ is satisfied and the Lagrange multipliers are known and fixed. Then
\[
\bar S=\hat S_1^{\pr\pr}(\lambda)=\hat S_2^{\pr\pr}(\lambda).
\]
\end{theorem}
\begin{proof}
It is obvious that $\hat S_1^{\pr\pr}(\lambda)\subseteq\hat S_2^{\pr\pr}(\lambda)$.

We prove that $S_2^{\pr\pr}(\lambda)\subset\bar S$. Let $x$ be an arbitrary point from the set $S_2^{\pr\pr}(\lambda)$. It follows from $\nabla g_i(\bar x)^T(x-\bar x)\ge 0$, $i\in\tilde I(\bar x,\lambda)$  and $(\ref{7})$ that $\nabla f(\bar x)^T(x-\bar x)\le 0$. Therefore $x\in \hat S_2^\pr(\lambda)$. We conclude from $\hat S_2^\pr(\lambda)=\bar S$ in Theorem \ref{th4} that $x\in\bar S$.

We prove that $\bar S\subseteq\hat S_1^{\pr\pr}(\lambda)$. Let $x$ be arbitrary point from $\bar S$. Due to the quasiconvexity of $f$ and $f(x)=f(\bar x)$ it follows from Lemma \ref{lema2} that $\nabla f(\bar x)^T(x-\bar x)\le 0$. Since $g_i$ are also quasiconvex and $g_i(x)\le 0=g_i(\bar x)$ for all $i\in I(\bar x)$, then by Lemma \ref{lema2}, we obtain that 
\[
\nabla g_i(x)^T(x-\bar x)\le 0 \quad\textrm{for all}\quad i\in I(\bar x).
\]
 By KKT conditions Equations (\ref{1}) are satisfied. Therefore, $\nabla g_i(x)^T(x-\bar x)=0$ for all $i\in\tilde I(\bar x,\lambda)$. Then, it follows from 
$\hat S_1^\pr(\lambda)=\bar S$ that $x\in \hat S_1^{\pr\pr}(\lambda)$.
\end{proof}

\section{Comparisons}
In this section, we compare our results with the previous ones. 

In Ref. \cite{suz15}, Suzuki and Kuroiwa derived characterizations of the solution set of a program with essentially quasiconvex objective function in terms of Greenberg-Pierskalla subdifferential.

A quasiconvex function is called essentially quasiconvex iff each local minimum is global. It is known that a continuous quasiconvex function is essentially quasiconvex if and only if it is semistrictly quasiconvex (see Theorem 3.37 in \cite{avr88}).

Greenberg-Pierskalla subdifferential \cite{gre73} for a quasiconvex function $f$ at the point $x_0$ is the following set:
\[
\partial^{GP} f(x_0)=\{v\in\R^n\mid v^T(x-x_0)\ge 0\quad\textrm{implies}\quad f(x)\ge f(x_0)\}.
\]
It is easy to see that the gradient of the function is not necessarily included in Greenberg-Pierskalla subdifferential and this subdifferential does not reduce to the gradient of the function, when the last one is continuously differentiable. For example, consider the continuously differentiable function of one variable defined as follows: $f(x)=x^2$, if $x\ge 0$, and $f(x)=-x^2$, if $x\le 0$. Then $\nabla f(0)=0$, but $\partial^{GP} f(0)=(0,+\infty)$. It follows from here that even the function is essentially quasiconvex, Theorem \ref{th1} does not follows from the results in \cite{suz15}. Moreover, we consider functions, which are not essentially quasiconvex.
%We continue the comparison with the characterizations in \cite{suz15} in the case when the function is differentiable. %We should explaim why the characterizations in \cite{suz15} are more usefull for problems, whose objective function is not differentiable. 
It is easy to calculate the derivative. We need to solve some system of equations to apply Theorem \ref{th1}. On the other hand, it is more difficult to construct the Greenberg-Pierskalla subdifferential than to calculate the derivative of the function, provided that the last one exists.

The following result is a part of Theorem 8 in \cite{suz15}:
\begin{proposition}\label{pr1}
Let $f$ be an upper semi-continuous essentially quasiconvex function, $S$ is a nonempty convex subset of $\R^n$, $\bar x\in\bar S$. Then
\[
\bar S=\{x\in S\mid\exists v\in\partial^{GP} f(\bar x)\cap\partial^{GP} f(x)\;\; \rm{such}\;\rm{that} \;\; v^T(x-\bar x)=0\}.
\]
\end{proposition}

Let us consider Example \ref{ex1} again. We can find the solution set using Greenberg-Pierskalla subdifferential.
Then $\partial^{GP} f(\bar x)=(0,+\infty)\times (0,+\infty)$.  Let $x=(x_1,x_2)$, where $x_1\ge 0$, $x_2\ge 0$, but $(x_1,x_2)\ne (0,0)$. Then 
\[
\partial^{GP} f(x)=\{y\in\R^2\mid  y_1=t x_1,\; y_2=t x_2, \; t>0\}.
\]
Let $x=(x_1,x_2)$, where $x_1>0$, $x_2<0$. Then
\[
\partial^{GP} f(x)=\{y\in\R^2\mid  y_1=t x_1,\; y_2=0, \; t>0\}.
\]
It is easy to see that in both cases the equality $v^T(x-\bar x)=0$ is not satisfied. Let $x=(x_1,x_2)$ where $x_1=0$, $x_2\le 0$.
Then $v^T(x-\bar x)=0$, where $v=(v_1,v_2)$, $v_1=t$, $v_2=0$, $t>0$. It follows from Proposition \ref{pr1} that the solution set is
$\bar S=\{x=(x_1,x_2)\mid x_1=0,\; x_2\le 0\}$.

In Ref. \cite{suz17}, Suzuki and Kuroiwa derived characterizations of the solution set of programs with non-essentially and essentially quasiconvex objective function in terms of Mar\-ti\-nez-Legaz subdifferential. Martinez-Legaz subdifferential \cite{mar88} of a function $f$ at a point $x\in\R^n$ is defined as follows:
\[
\partial^M f(x)=\{(v,t)\in\R^{n+1}\mid \inf_{v^T y\ge t}f(y)\ge f(x),\; v^T x\ge t\}.
\]
They derived the following result:
\begin{proposition}\label{suzuki}
Let $f$ be a usc quasiconvex function and $\bar x$ be a global solution of the problem (P). Then the following sets are equal:

\begin{tabular}{@{}rl}
(i) & $\bar S=\{x\in S\mid f(x)=\min_y\in S\}$, \\
(ii) & $M_1=\{x\in S\mid\partial^M f(\bar x)\cap\partial^M f(x)\ne\emptyset\}$, \\
(iii) & $M_2=\{x\in S\mid\exists (v,t)\in\partial^M f(x)\;\;$ such that $\;\; v^T\bar x\ge t\}$.
\end{tabular}
%\begin{enumerate}
%\item[(i)]
%$\bar S=\{x\in S\mid f(x)=\min_{y\in S} f(y)\}$,
%\item[(ii)]
%$M_1=\{x\in S\mid\partial^M f(\bar x)\cap\partial^M f(x)\ne\emptyset\}$,
%\item[(iii)]
%$M_2=\{x\in S\mid\exists (v,t)\in\partial^M f(x)\;\; \rm{such}\;\rm{that} \;\; v^T\bar x\ge t\}$.
%\end{enumerate}
\end{proposition} 

It is not a very easy task to compute Martinez-Legaz subdifferential for functions which depend on two or more variables. If we want to apply these characterizations, then we should compute Martinez-Legaz subdifferential at every point $x\in S$. To construct this subdifferential we should solve a minimization problem. The easiest case is the one-dimensional problem with a monotone function, because in this case we can easy compute the infimum. Consider the following example, which is similar to the Example 3.1 in \cite{suz17}. Really, this example is different from \cite[Example 3.1]{suz17}, because the function in Example 3.1 is not continuously differentiable.

\begin{example}
Let $S=[0,2]$ and $f$ be a real-valued function on $R$ as follows:
\[
f(x)=\left\{\begin{array}{rl}
-x^2, & x\in(-\infty,0),\\
0, & x\in[0,1],\\
(x-1)^2, & x\in (1,+\infty)
\end{array}\right.
\]
Let $\bar x=0$ be a known solution. The solution set is the interval $[0,1]$. We can check this very easy by Theorem \ref{th2}, because $\nabla f(x)=0$ if $x\in [0,1]$ and $\nabla f(x)=2(x-1)\ne 0$ is $x\in (1,2]$. The function is pseudoconvex at every point $x\in (1,2]$.

Let us solve the problem applying Theorem \ref{suzuki}. First, we compute Martinez-Legaz subdifferential when $x\in [0,1]$. Let $(v,t)\in\partial^M f(x)$. Therefore $\inf\{f(y)\mid v y\ge t\}\ge f(x)=0$.  We prove that $v>0$. Arguing by contradiction, if $v<0$, then $\inf\{f(y)\mid v y\ge t\}=-\infty$, which is a contradiction. If $v=0$ and $t\le 0$, then $\inf\{f(y)\mid v y\ge t\}=-\infty$, a contradiction again. If $v=0$ and $t>0$, then $0=v x\ge t>0$. Impossible. Hence $v>0$. By the definition of Martinez-Legaz subdifferential  we have
\[
f(x)\ge\inf\{f(y)\mid v y\ge t\}=\inf\{f(y)\mid y\ge t/v\}\ge f(x)=0.
\]
It follows from here that 
\[
\inf\{f(y)\mid v y\ge t\}=\inf\{f(y)\mid y\ge t/v\}=0
\]
and $0\le t/v\le 1$. Thus
\[
\partial^M f(x)\subseteq \{(v,t)\mid v>0,\; 0\le t\le v,\; v x\ge t\}.
\]
We prove the inverse inclusion. Let $v>0$, $0\le t\le v$, and $v x\ge t$. Therefore
\[
\inf\{f(y)\mid v y\ge t\}=\inf\{f(y)\mid y\ge t/v\}=0=f(x).
\]
We compute the Martinez-Legaz subdifferential when $x\in (1,2]$. We prove that $v>0$. If $v<0$ or $v=0$, $t\le 0$, then $\inf\{f(y)\mid v y\ge t\}=-\infty$. If $v=0$ and $t>0$, then $0=v x\ge t>0$. Impossible. Hence $v>0$. We have
\[
f(x)\ge\inf\{f(y)\mid v y\ge t\}=\inf\{f(y)\mid y\ge t/v\}\ge f(x).
\]
Therefore $v x=t$ and \[
\partial^M f(x)\subseteq \{(v,vx)\mid v>0\}.
\]
We prove the inverse inclusion. Let $v>0$. Therefore
\[
\inf\{f(y)\mid v y\ge t\}=\inf\{f(y)\mid y\ge t/v\}=f(x).
\]
At last, we find the solution set. If there exists $(v,t)\in\partial^M f(\bar x)\cap\partial^M f(x)$, such that $1<x\le 2$, then by $vx=t$ and $0\le t\le v$ we conclude that $x\le 1$, a contradiction. Therefore, by Proposition \ref{suzuki}, $\bar S=[0,1]$.
\end{example}

Consider the sets:
\[\begin{array}{l}
\hat T_1:=\{x\in S\mid\nabla f(\bar x)^T(x-\bar x)=0,\; \exists p(x)>0: \nabla f(x)=p(x)\nabla f(\bar x)\}; \\
%\]
%\[
\hat T_2:=\{x\in S\mid\nabla f(\bar x)^T(x-\bar x)\le 0,\; \exists p(x)>0: \nabla f(x)=p(x)\nabla f(\bar x)\}.
\end{array}\]

The following result is a corollary of Theorem \ref{th1}. It is a particular case of Theorem \ref{th1} by Ivanov in \cite{jota2013}, where the objective function is differentiable and pseudoconvex:

\begin{corollary}\label{cor1}
Let $\Gamma\subseteq\E$ be an open convex set, $S\subseteq\Gamma$ be a convex one and $f:\Gamma\to\R$ be a continuously differentiable pseudoconvex function. Suppose that $\bar x\in\bar S$ is a known solution of {\rm (P)}. Then
\[
\bar S=\hat T_1=\hat T_2.
\]
\end{corollary}
 \begin{proof}  
If $\nabla f(\bar x)\ne 0$, then it follows from Theorem \ref{th1} that $\bar S=\hat S_1=\hat S_2$.  Suppose that there exists $x\in S$ with
\[
\nabla f(\bar x)^T(x-\bar x)\le 0,\; \exists p>0: \nabla f(x)=p\nabla f(\bar x),\;\nabla f(x)=0.
\]
It follows from here that $\nabla f(\bar x)=0$, which is a contradiction. Therefore $\bar S=\hat T_1=\hat T_2$.

Let $\nabla f(\bar x)=0$. Then $\hat T_1=\hat T_2=\{x\in S\mid\nabla f(x)=0\}$. By pseudoconvexity $\bar x$ is a global minimizer of $f$ on $\Gamma$ and
\[
\bar S=\{x\in S\mid f(x)=f(\bar x)\}=\{x\in S\mid\nabla f(x)=0\}.
\]
The proof is complete.
\end{proof}

\begin{corollary}
Suppose additionally to the hypothesis of Theorem \ref{th1} that $f$ is convex. Then, $\nabla f(x)=\nabla f(\bar x)$, provided that $x\in\bar S$.
\end{corollary}
\begin{proof}
It follows from convexity that $f(y)-f(x)\ge\nabla f(x)^T(y-x)$ for all $y\in S$. According to the relation $\bar S=\hat S_1$, we obtain that $\nabla f(x)=p\nabla f(\bar x)$, where $p=\norm{\nabla f(x)}/\norm{\nabla f(\bar x)}$ and $\nabla f(\bar x)^T(x-\bar x)=0$. Also $f(x)=f(\bar x)$ by the inclusion $x\in \bar S$. Therefore, the following inequality holds:
\[
f(y)-f(\bar x)\ge p\nabla f(\bar x)^T(y-\bar x+\bar x-x)=p\nabla f(\bar x)^T(y-\bar x).
\]
It follows from here that $p\nabla f(\bar x)$ belongs to the subdifferential of $f$ at $\bar x$, denoted by $\partial f(\bar x)$. On the other hand, by continuous differentiability, the subdifferential contains the unique element $\nabla f(\bar x)$. Hence, $p=1$ and $\nabla f(x)=\nabla f(\bar x)$.
\end{proof}

It follows from this corollary that Theorem \ref{th1} together with Theorem \ref{th2} are  generalization of Theorem 1 in the paper by Mangasarian \cite{man88}, where the function is twice continuously differentiable and convex.

Consider the following sets:
\[\begin{array}{l}
T_1:=\{x\in S\mid \nabla f(x)^T(\bar x-x)=0\}, \\
%\]
%\[
T_2:=\{x\in S\mid \nabla f(x)^T(\bar x-x)\ge 0\}, \\
%\]
%\[
T_3:=\{x\in S\mid \nabla f(x)^T(\bar x-x)=\nabla f(\bar x)^T(x-\bar x)\}, \\
%\]
%\[
T_4:=\{x\in S\mid \nabla f(x)^T(\bar x-x)\ge\nabla f(\bar x)^T(x-\bar x)\}, \\
%\]
%\[
T_5:=\{x\in S\mid\nabla f(x)^T(\bar x-x)=\nabla f(\bar x)^T(x-\bar x)=0\}.
\end{array}\]

The following result is a corollary of Theorem \ref{th3}. It is a particular case of \cite[Theorem 4.1]{2001}, \cite[Theorem 3.1]{2003}:

\begin{corollary}\label{cor2}
 Let $\Gamma\subseteq\E$ be an open convex set, $S\subseteq\Gamma$ be an arbitrary convex one and the function $f$ be defined on $\Gamma$.
Suppose $\bar x\in\bar S$ is a known solution of {\rm (P)}.
If the function $f$ is continuously differentiable and  pseudoconvex on $\Gamma$, then
\[
\bar  S=T_1=T_2=T_3=T_4=T_5.
\]
\end{corollary}
\begin{proof} 
It is obvious that $S_i\subseteq T_i$, $i=1,2,3,4,5$ and $T_5\subseteq T_1\subseteq T_2$ and $T_5\subseteq T_3\subseteq T_4$.

The proof of $T_4\subseteq T_2$ uses the arguments of Theorem \ref{th3}.

We prove that  $\bar  S\subseteq T_5$. If $\nabla f(\bar x)\ne 0$, then by Theorem \ref{th3}, we obtain that $\bar S=S_5\subseteq T_5$. Let $\nabla f(\bar x)=0$. Then by pseudoconvexity $\bar x$ is global minimizer of $f$ over $\Gamma$. It follows from here that $\bar S=\{x\in S\mid\nabla f(x)=0\}\subseteq T_5$.

At last, we prove that $T_2\subseteq\bar S$. Let $x\in T_2$. Assume the contrary that $x\notin\bar S$. It follows from here that 
$f(\bar x)<f(x)$. By pseudoconvexity, we obtain that $\nabla f(x)^T(\bar x-x)<0$, which is a contradiction to the assumption $x\in T_2$.
\end{proof}

%Consider the following result due to Jeyakumar and Yang \cite[Proposition 3.1]{jey95}:
%\begin{proposition}
%Let $f$ be continuously differentiable and pseudolinear on an open convex set containing the convex set $S\subseteq\R^n$; let 
%$\bar x\in\bar S$. Suppose that $\nabla f(x)\ne 0$ on $S$. Then
%\[
%\bar S\subseteq\{x\in S\mid\nabla f(x)/\norm{\nabla f(x)}=\nabla f(\bar x)/\norm{\nabla f(\bar x)}\}\subseteq\{x\in S\mid\nabla f(x)^T(\bar x-x)\le 0\}.
%\]
%\end{proposition}

%Theorem \ref{th3} is a generalization of the main results  in the paper \cite{jey95} by Jeyakumar and Yang (see Theorem 3.1, Theorem 3.2 and Corollary 3.1), where the objective function is differentiable and pseudolinear.

\end{document}